
\baselineskip=14pt
\parskip=10pt

\magnification=\magstephalf

\def\1{{\overline{1}}}
\def\2{{\overline{2}}}
\parindent=0pt
\overfullrule=0in

\def\frac#1#2{{#1 \over #2}}
\centerline
{\bf The Average AMS fellow who died before Aug. 2022 lived 83.36 years,} 
\centerline
{\bf published 98.53 papers, and was cited 2473.23 times}
\bigskip
\centerline
{\it Livia A. STONE and Doron ZEILBERGER}

\bigskip

{\bf Abstract}:  Using data collected from the AMS websites, and the internet in general, we compute the weight-enumerator of
the set of $145$ AMS fellows who passed away before Aug. 2022, according to the statistics [Age, Number of Publications, Number of Citations].
Using Maple and elementary calculus, we deduced from it the numbers stated in the title, as well as much more
detailed and fascinating statistical information. For example, while the correlation coefficient between the number of
publications and number of citations is, not surprisingly, positive, it is far lower than one would expect, being a mere $0.5797$ .
The appendices contain the rankings according to longevity, number of publications, and number of citations.

{\bf Lemma}: The weight-enumerator of the set of AMS fellows who died before Aug. 2022, 
where the weight of a fellow is
$$
a^{AgeAtDeath}\,p^{NumberOfPublications}\,c^{NumberOfCitations} \quad,
$$
let's call it $F(a,p,c)$, is the following polynomial of degree $100$ in $a$, degree $453$ in $p$, and degree $27492$ in $c$:

$$
F(a,p,c)=\,
$$
$$
{c}^{173}{p}^{24}{a}^{48}+ \left( {c}^{2898}{p}^{122}+{c}^{1859}{p}^{111} \right) {a}^{58}+ \left( {c}^{2116}{p}^{48}+{c}^{1265}{p}^{57} \right) {a}^{59}+{c}^{8482}{p}^{113}{a}^{60}+{c}^{101}{p}^{40}{a}^{63}
$$
$$
+ \left( {c}^{9810}{p}^{453}+{c}^{882}{p}^{46} \right) {a}^{65}+ \left( {c}^{2084}{p}^{84}+{c}^{1337}{p}^{92}+{c}^{366}{p}^{39} \right) {a}^{66}+{c}^{3924}{p}^{117}{a}^{68}+{c}^{756}{p}^{55}{a}^{70}
$$
$$
+{c}^{11987}{p}^{323}{a}^{72}+ \left( {c}^{1692}{p}^{133}+{c}^{734}{p}^{50} \right) {a}^{73}+ \left( {c}^{6030}{p}^{122}+{c}^{2677}{p}^{186}+{c}^{1624}{p}^{77}+{c}^{341}{p}^{18} \right) {a}^{74}
$$
$$
+ \left( {c}^{749}{p}^{47}+{c}^{658}{p}^{42}+{c}^{145}{p}^{25} \right) {a}^{75}+ \left( {c}^{3918}{p}^{106}+{c}^{1867}{p}^{107}+{c}^{1781}{p}^{136}+{c}^{1309}{p}^{183} \right) {a}^{76}
$$
$$
+ \left( {c}^{2137}{p}^{175}+{c}^{451}{p}^{45} \right) {a}^{77}+ \left( {c}^{4633}{p}^{189}+{c}^{4068}{p}^{132}+{c}^{1891}{p}^{66}+{c}^{914}{p}^{61} \right) {a}^{78}
$$
$$
+ \left( {c}^{2563}{p}^{48}+{c}^{1596}{p}^{65}+{c}^{1550}{p}^{70}+{c}^{181}{p}^{34}+{c}^{157}{p}^{29} \right) {a}^{79}
$$
$$
+ \left( {c}^{3732}{p}^{168}+{c}^{1370}{p}^{49}+{c}^{503}{p}^{33}+{c}^{428}{p}^{23}+{c}^{281}{p}^{37} \right) {a}^{80}
$$
$$
+ \left( {c}^{15121}{p}^{159}+{c}^{1937}{p}^{166}+{c}^{873}{p}^{31}+{c}^{560}{p}^{67}+{c}^{121}{p}^{39} \right) {a}^{81}
$$
$$
+ \left( {c}^{4096}{p}^{181}+{c}^{2776}{p}^{138}+{c}^{2336}{p}^{149}+{c}^{2111}{p}^{190}+{c}^{2035}{p}^{64}+{c}^{1507}{p}^{84}+{c}^{1281}{p}^{164} \right) {a}^{82}
$$
$$
+ \left( {c}^{2604}{p}^{79}+{c}^{1746}{p}^{153}+{c}^{1658}{p}^{81}+{c}^{1509}{p}^{78}+{c}^{781}{p}^{26}+{c}^{133}{p}^{27} \right) {a}^{83}
$$
$$
+ \left( {c}^{3509}{p}^{118}+{c}^{3484}{p}^{82}+{c}^{1497}{p}^{146}+{c}^{1510}{p}^{53}+{c}^{1172}{p}^{74}+{c}^{542}{p}^{46}+{c}^{526}{p}^{34}+{c}^{336}{p}^{23}+{c}^{161}{p}^{31} \right) {a}^{84}
$$
$$
+ \left( {c}^{7500}{p}^{363}+{c}^{4772}{p}^{193}+{c}^{4490}{p}^{126}+{c}^{2116}{p}^{164}+{c}^{1231}{p}^{47} \right) {a}^{85}
$$
$$
+ \left( {c}^{4802}{p}^{163}+{c}^{2165}{p}^{57}+{c}^{1826}{p}^{66}+{c}^{1118}{p}^{71}+{c}^{452}{p}^{46}+{c}^{275}{p}^{42} \right) {a}^{86}
$$
$$
+ \left( {c}^{27492}{p}^{237}+{c}^{3185}{p}^{29}+{c}^{2519}{p}^{71}+{c}^{1642}{p}^{105}+{c}^{1211}{p}^{31}+{c}^{978}{p}^{41} \right) {a}^{87}
$$
$$
+ \left( {c}^{2845}{p}^{102}+{c}^{1850}{p}^{106}+{c}^{1726}{p}^{69}+{c}^{1173}{p}^{49}+{c}^{1146}{p}^{29}+{c}^{1086}{p}^{81}+{c}^{426}{p}^{35}+{c}^{142}{p}^{19}+{c}^{12}{p}^{11} \right) {a}^{88}+ 
$$
$$
( {c}^{6664}{p}^{115}+{c}^{5937}{p}^{180}+{c}^{4743}{p}^{200}+{c}^{4367}{p}^{114}+{c}^{4066}{p}^{167}+{c}^{3814}{p}^{271}+{c}^{2248}{p}^{101}+{c}^{1610}{p}^{89}+{c}^{1130}{p}^{159}
$$
$$
+{c}^{1107}{p}^{111}+{c}^{1054}{p}^{61}+{c}^{791}{p}^{87}+
{c}^{239}{p}^{99} ) {a}^{89}+ 
$$
$$
\left( {c}^{12281}{p}^{269}+{c}^{3685}{p}^{110}+{c}^{2851}{p}^{222}+{c}^{145}{p}^{37} \right) {a}^{90}
$$
$$
+ \left( {c}^{2639}{p}^{288}+{c}^{1183}{p}^{57}+{c}^{664}{p}^{118}+{c}^{525}{p}^{39}+{c}^{152}{p}^{17}+{c}^{52}{p}^{36} \right) {a}^{91}
$$
$$
+ \left( {c}^{1601}{p}^{59}+{c}^{1065}{p}^{59}+{c}^{498}{p}^{99} \right) {a}^{92}
$$
$$
+ \left( {c}^{8722}{p}^{114}+{c}^{4107}{p}^{70}+{c}^{3655}{p}^{321}+{c}^{3633}{p}^{105}+{c}^{2886}{p}^{87}+{c}^{2128}{p}^{178}+{c}^{1254}{p}^{37}+{c}^{41}{p}^{19} \right) {a}^{93}
$$
$$
+ \left( {c}^{5626}{p}^{95}+{c}^{1295}{p}^{82}+{c}^{214}{p}^{27}+{c}^{2}{p}^{5} \right) {a}^{94}+ \left( {c}^{17621}{p}^{175}+{c}^{586}{p}^{58}+{c}^{400}{p}^{182} \right) {a}^{95}
$$
$$
+ \left( {c}^{1632}{p}^{90}+{c}^{712}{p}^{46}+{c}^{206}{p}^{66} \right) {a}^{96}+
$$
$$
 \left( {c}^{5215}{p}^{96}+{c}^{778}{p}^{75} \right) {a}^{97}+{c}^{272}{p}^{86}{a}^{99}+{c}^{401}{p}^{93}{a}^{100} \quad .
$$

First we turn the enumeration generating function to a probability generating function by dividing by $F(1,1,1)=145$.

$$
f (a,b,c):=\frac{F(a,b,c)}{F(1,1,1)}= \frac{1}{145} F(a,b,c) \quad .
$$

We are now ready to prove our theorem.

{\bf Theorem }: 

$\bullet$ The average age at death of an AMS Fellow who died before Aug. 2022 is  $83.3655$,
the standard deviation is $9.343$, the skewness is, $-1.079$, and the kurtosis is $4.3277$.

$\bullet$ 
The average number of publications is $98.531$,
the standard deviation is $73.942$, the skewness is $1.757$, and the kurtosis is $7.1779$.

$\bullet$
The average number of citations is  $2473.234$,
the standard deviation is $3455.194$,  the skewness is $3.97$,
and the kurtosis is, $23.993$.

$\bullet$ The correlation coefficient between the age at death and the number of publications is $0.003$,
between the age at death and the number of citations is $0.011$, while the
correlation between the number of publications and number of citations is $0.57973$.

{\bf Proof}: 
We have, of course
$$
\mu_{age}= \frac{\partial}{\partial a} f(a,p,c)\bigl \vert_{a=1,p=1.c=1} \quad,
$$

$$
\mu_{pub}= \frac{\partial}{\partial p} f(a,p,c)\bigl \vert_{a=1,p=1.c=1} \quad,
$$

$$
\mu_{cit}= \frac{\partial}{\partial c} f(a,p,c)\bigl \vert_{a=1,p=1.c=1} \quad .
$$

Now, let's define
$$
{\bar f}(a,p,c) := \frac{f(a,p,c)}{a^{\mu_{age}} \, p^{\mu_{pub}} \, c^{\mu_{cit}}} \quad.
$$

We have, of course, that the {\bf moments about the mean} are:
$$
M_{r,s,t}=E[(age-\mu_{age})^r\,(pub-\mu_{pub})^s\,(cit-\mu_{cit})^t]=  \left (a \frac{\partial}{\partial a} \right)^r
\left (p \frac{\partial}{\partial p} \right)^s
\left (c \frac{\partial}{\partial c} \right)^t {\bar f}(p,a,c) \Bigl \vert_{a=1,p=1,c=1} \quad ,
$$
and finally the {\bf scaled moments} are
$$
\alpha_{r,s,t} \, = \, \frac{M_{r,s,t}}{M_{2,0,0}^{r/2}M_{0,2,0}^{s/2}M_{0,0,2}^{t/2}} \quad.
$$

The numbers stated in the theorem are, respectively, 

$\bullet$ $\mu_{age}$, $\sqrt{M_{2,0,0}}$, $\alpha_{3,0,0}$, $\alpha_{4,0,0}$

$\bullet$ $\mu_{pub}$, $\sqrt{M_{0,2,0}}$, $\alpha_{0,3,0}$, $\alpha_{0,4,0}$

$\bullet$ $\mu_{cit}$, $\sqrt{M_{0,0,2}}$, $\alpha_{0,0,3}$, $\alpha_{0,0,4}$

$\bullet$ $\alpha_{1,1,0}$,  $\alpha_{1,0,1}$,  $\alpha_{0,1,1}$ .

{\bf Appendix I: The Dead AMS Fellows According to Age}

Here are all the dead fellows listed, in decreasing order of their age at 
death, followed by their number of publications and number of citations in MathSciNet (viewed Aug. 2022).

               There was 1 fellow who died at the age of {\bf 100}

                       Maurice H. Heins ( 1915-2015) 93;401

               There was 1 fellow who died at the age of {\bf 99}

                         Lee Lorch ( 1915-2014) 86;272

               There were 2 fellows who died at the age of {\bf 97}

                 Komaravolu Chandrasekharan ( 1920-2017) 75;778

                      Isadore M. Singer ( 1924-2021) 96;5215

               There were 3 fellows who died at the age of {\bf 96}

                     Jane Cronin Scanlon ( 1922-2018) 66;206

                      Bjarni Jonsson ( 1920-2016) 90;1632

                      Richard D. Schafer ( 1918-2014) 46;712

               There were 3 fellows who died at the age of {\bf 95}

                        Edgar H. Brown ( 1926-2021) 58;586

                      Donald A.S. Fraser ( 1925-2020) 182;400

                     Louis Nirenberg ( 1925-2020) 175;17621

               There were 4 fellows who died at the age of {\bf 94}

                         Lida K. Barrett ( 1927-2021) 5;2

                  Cathleen Synge Morawetz ( 1923-2017) 82;1295

                        John T. Tate ( 1925-2019) 95;5626

                        John H. Walter ( 1927-2021) 27;214

               There were 8 fellows who died at the age of {\bf 93}

                       Tom M. Apostol ( 1923-2016) 87;2886

                   Arthur Herbert Copeland ( 1926-2019) 19;41

                     Richard V. Kadison ( 1925-2018) 105;3633

                      Joseph B. Keller ( 1923-2016) 321;3655

                        John C. Moore ( 1923-2016) 37;1254

                    Murray Rosenblatt ( 1926-2019) 178;2128

                   Lloyd Stowell Shapley ( 1923-2016) 70;4107

                       Guido Weiss ( 1928-2021) 114;8722

               There were 3 fellows who died at the age of {\bf 92}

                       Edward Fadell ( 1926-2018) 59;1065

                        Walter Noll ( 1925-2017) 59;1601

                        Herman Rubin ( 1926-2018) 99;498

               There were 6 fellows who died at the age of {\bf 91}

                        Harvey Cohn ( 1923-2014) 118;664

                    Jean-Pierre Kahane ( 1926-2017) 288;2639

                      Arthur P. Mattuck ( 1930-2021) 17;152

                       Isaac Namioka ( 1928-2019) 57;1183

                        Howard Osborn ( 1928-2019) 36;52

                        Maurice Sion ( 1927-2018) 39;525

               There were 4 fellows who died at the age of {\bf 90}

                      Steve Armentrout ( 1930-2020) 37;145

                    Sir Michael Atiyah ( 1929-2019) 269;12281

                    Manfredo P. doCarmo ( 1928-2018) 110;3685

                      Eugene B. Dynkin ( 1924-2014) 222;2851

               There were 13 fellows who died at the age of {\bf 89}

                      Felix E. Browder ( 1927-2016) 200;4743

                      Jim Douglas, Jr. ( 1927-2016) 180;5937

                     Branko Grunbaum ( 1929-2018) 271;3814

                       Junichi Igusa ( 1924-2013) 89;1610

                     Bertram Kostant ( 1928-2017) 114;4367

                   Wilhelmus Luxemburg ( 1929-2018) 111;1107

                      Sibe Mardešic ( 1927-2016) 159;1130

                     Albert Nijenhuis ( 1926-2015) 61;1054

                        Vera Pless ( 1931-2020) 101;2248

                      Mary Ellen Rudin ( 1924-2013) 87;791

                      Victor L. Shapiro ( 1924-2013) 99;239

                     Hans Weinberger ( 1928-2017) 115;6664

                       Harold Widom ( 1932-2021) 167;4066

               There were 9 fellows who died at the age of {\bf 88}

                        William Bade ( 1924-2012) 35;426

                      Richard L. Bishop ( 1931-2019) 49;1173

                     Solomon Feferman ( 1928-2016) 29;1146

                      Ronald K. Getoor ( 1929-2017) 106;1850

                     Heini Halberstam ( 1926-2014) 81;1086

                       James M. Kister ( 1930-2018) 19;142

                     Carole B. Lacampagne ( 1933-2021) 11;12

                     Ernest A. Michael ( 1925-2013) 102;2845

             Conjeeveram Srirangachari Seshadri ( 1932-2020) 69;1726

               There were 6 fellows who died at the age of {\bf 87}

                   Elliott Ward Cheney, Jr. ( 1929-2016) 105;1642

                        Robert Ellis ( 1926-2013) 41;978

                        John Hempel ( 1935-2022) 31;1211

                       John F. Nash, Jr. ( 1928-2015) 29;3185

                      Robert R. Phelps ( 1926-2013) 71;2519

                      Elias M. Stein ( 1931-2018) 237;27492

               There were 6 fellows who died at the age of {\bf 86}

                      Richard Askey ( 1933-2019) 163;4802

                    Donald L. Burkholder ( 1927-2013) 57;2165

                     Pierre E. Conner, Jr. ( 1932-2018) 71;1118

                        Louis Howard ( 1929-2015) 42;275

                       Rudolf Kalman ( 1930-2016) 66;1826

                         Ray Kunze ( 1928-2014) 46;452

               There were 5 fellows who died at the age of {\bf 85}

                         Roy Adler ( 1931-2016) 47;1231

                       Peter Duren ( 1935-2020) 126;4490

                      Ronald Graham ( 1935-2020) 363;7500

                      Andras Hajnal ( 1931-2016) 164;2116

                        Harry Kesten ( 1931-2016) 193;4772

               There were 9 fellows who died at the age of {\bf 84}

                       Donald Babbitt ( 1936-2020) 31;161

                      Robert J. Blattner ( 1931-2015) 23;336

                      Harold M. Edwards ( 1936-2020) 46;542

                      Edward G. Effros ( 1935-2019) 82;3484

                        Paul Fife ( 1930-2014) 118;3509

                      Solomon Golomb ( 1932-2016) 146;1497

                  Jacques Jean-Pierre Neveu ( 1932-2016) 74;1172

                      Robert T. Seeley ( 1932-2016) 53;1510

                        Henry Wente ( 1936-2020) 34;526

               There were 6 fellows who died at the age of {\bf 83}

                     JMichael Boardman ( 1938-2021) 26;781

                     W. Wistar Comfort ( 1933-2016) 153;1746

                      Richard Pollack ( 1935-2018) 81;1658

                      Charles J. Stone ( 1936-2019) 79;2604

                       Joseph F. Traub ( 1932-2015) 78;1509

                     Alphonse T. Vasquez ( 1938-2021) 27;133

               There were 7 fellows who died at the age of {\bf 82}

                    Shreeram Abhyankar ( 1930-2012) 190;2111

                     Richard M. Dudley ( 1938-2020) 138;2776

                   Clifford John Earle, Jr. ( 1935-2017) 84;1507

                      Mark Mahowald ( 1931-2013) 164;1281

                       Edward Nelson ( 1932-2014) 64;2035

                      Joel A. Smoller ( 1935-2017) 181;4096

                 Veeravalli S. Varadarajan ( 1937-2019) 149;2336

               There were 5 fellows who died at the age of {\bf 81}

                       Kenneth I. Appel ( 1932-2013) 31;873

                     Gilbert Baumslag ( 1933-2014) 166;1937

                       Samuel Gitler ( 1933-2014) 67;560

                     Lars Hormander ( 1931-2012) 159;15121

                    Aderemi Oluyomi Kuku ( 1941-2022) 39;121

               There were 5 fellows who died at the age of {\bf 80}

                     Ronald G. Douglas ( 1938-2018) 168;3732

                       Robert M. Miura ( 1938-2018) 49;1370

                       Paul J. Sally, Jr. ( 1933-2013) 37;281

                       Charles C. Sims ( 1937-2017) 33;503

                       Myles Tierney ( 1937-2017) 23;428

               There were 5 fellows who died at the age of {\bf 79}

                        Alan Baker ( 1939-2018) 65;1596

                      Elwyn Berlekamp ( 1940-2019) 70;1550

                    Aldridge Bousfield ( 1941-2020) 48;2563

                       Kenneth I. Gross ( 1938-2017) 34;181

                       Lesley M. Sibner ( 1934-2013) 29;157

               There were 4 fellows who died at the age of {\bf 78}

                     Joyce R. McLaughlin ( 1939-2017) 61;914

                       George R. Sell ( 1937-2015) 132;4068

                    Robert Strichartz ( 1943-2021) 189;4633

                      William A. Veech ( 1938-2016) 66;1891

               There were 2 fellows who died at the age of {\bf 77}

                        David Cantor ( 1935-2012) 45;451

                      Lawrence Shepp ( 1936-2013) 175;2137

               There were 4 fellows who died at the age of {\bf 76}

                     Frederick R. Cohen ( 1946-2022) 136;1781

                    Vladimir F. Demyanov ( 1938-2014) 183;1309

                      Peter M. Gruber ( 1941-2017) 107;1867

                     Thomas M. Liggett ( 1944-2020) 106;3918

               There were 3 fellows who died at the age of {\bf 75}

                       Paul Chernoff ( 1942-2017) 42;658

                     Samuel M. Rankin, III ( 1945-2020) 25;145

                   Clarence W. Wilkerson, Jr. ( 1944-2019) 47;749

               There were 4 fellows who died at the age of {\bf 74}

                      Colin J. Bushnell ( 1947-2021) 77;1624

                      Anatole Katok ( 1944-2018) 122;6030

                      Lynn Arthur Steen ( 1941-2015) 18;341

                       Jan C. Willems ( 1939-2013) 186;2677

               There were 2 fellows who died at the age of {\bf 73}

                     Georgia Benkart ( 1949-2022) 133;1692

                     Oscar E. Lanford, III ( 1940-2013) 50;734

               There was 1 fellow who died at the age of {\bf 72}

                      Andrew J. Majda ( 1949-2021) 323;11987

               There was 1 fellows who died at the age of {\bf 70}

                       Steven M. Zucker ( 1949-2019) 55;756

               There was 1 fellow who died at the age of {\bf 68}

                     Vaughan F. R. Jones ( 1952-2020) 117;3924

               There were 3 fellows who died at the age of {\bf 66}

                       Walter Craig ( 1953-2019) 84;2084

                     Stephen A. Mitchell ( 1951-2017) 39;366

                       Edward Odell ( 1947-2013) 92;1337

               There were 2 fellows who died at the age of {\bf 65}

                 Jonathan Michael Borwein ( 1951-2016) 453;9810

                         David Goss ( 1952-2017) 46;882

               There was 1 fellows who died at the age of {\bf 63}

                       Lynne H. Walling ( 1958-2021) 40;101

               There was 1 fellow who died at the age of {\bf 60}

                    Andrei V. Zelevinsky ( 1953-2013) 113;8482

               There were 2 fellows who died at the age of {\bf 59}

                       Tim D. Cochran ( 1955-2014) 57;1265

                         John Roe ( 1959-2018) 48;2116

               There were, 2, fellows who died at the age of {\bf 58}

                     Nikolai Chernov ( 1956-2014) 111;1859

                       Robin Thomas ( 1962-2020) 122;2898

               There was 1 fellow who died at the age of {\bf 48}

                       Thomas Nevins ( 1972-2020) 24;173

{\bf Appendix II: Ranking the  AMS Fellows who died before Aug. 2022 According to Decreasing Number of Publications}

Here is the list of all 145 dead AMS fellows in decreasing order of their number of publications. 
It lists the year of birth, year of death, number of publications, and number of citations.

               1.  Jonathan Michael Borwein ( 1951-2016) 453;9810

                    2.  Ronald Graham ( 1935-2020) 363;7500

                    3.  Andrew J. Majda ( 1949-2021) 323;11987

                    4.  Joseph B. Keller ( 1923-2016) 321;3655

                  5.  Jean-Pierre Kahane ( 1926-2017) 288;2639

                   6.  Branko Grunbaum ( 1929-2018) 271;3814

                  7.  Sir Michael Atiyah ( 1929-2019) 269;12281

                    8.  Elias M. Stein ( 1931-2018) 237;27492

                    9.  Eugene B. Dynkin ( 1924-2014) 222;2851

                   10.  Felix E. Browder ( 1927-2016) 200;4743

                    11.  Harry Kesten ( 1931-2016) 193;4772

                 12.  Shreeram Abhyankar ( 1930-2012) 190;2111

                  13.  Robert Strichartz ( 1943-2021) 189;4633

                    14.  Jan C. Willems ( 1939-2013) 186;2677

                 15.  Vladimir F. Demyanov ( 1938-2014) 183;1309

                   16.  Donald A.S. Fraser ( 1925-2020) 182;400

                    17.  Joel A. Smoller ( 1935-2017) 181;4096

                    18.  Jim  Douglas, Jr. ( 1927-2016) 180;5937

                  19.  Murray Rosenblatt ( 1926-2019) 178;2128

                   20.  Lawrence Shepp ( 1936-2013) 175;2137

                  21.  Louis Nirenberg ( 1925-2020) 175;17621

                   22.  Ronald G. Douglas ( 1938-2018) 168;3732

                    23.  Harold Widom ( 1932-2021) 167;4066

                  24.  Gilbert Baumslag ( 1933-2014) 166;1937

                    25.  Mark Mahowald ( 1931-2013) 164;1281

                   26.  Andras Hajnal ( 1931-2016) 164;2116

                    27.  Richard Askey ( 1933-2019) 163;4802

                   28.  Sibe Mardesic ( 1927-2016) 159;1130

                  29.  Lars Hormander ( 1931-2012) 159;15121

                   30.  W. Wistar Comfort ( 1933-2016) 153;1746

               31.  Veeravalli S. Varadarajan ( 1937-2019) 149;2336

                   32.  Solomon Golomb ( 1932-2016) 146;1497

                   33.  Richard M. Dudley ( 1938-2020) 138;2776

                  34.  Frederick R. Cohen ( 1946-2022) 136;1781

                   35.  Georgia Benkart ( 1949-2022) 133;1692

                    36.  George R. Sell ( 1937-2015) 132;4068

                     37.  Peter Duren ( 1935-2020) 126;4490

                    38.  Robin Thomas ( 1962-2020) 122;2898

                    39.  Anatole Katok ( 1944-2018) 122;6030

                      40.  Paul Fife ( 1930-2014) 118;3509

                     41.  Harvey Cohn ( 1923-2014) 118;664

                   42.  Vaughan F. R. Jones ( 1952-2020) 117;3924

                   43.  Hans Weinberger ( 1928-2017) 115;6664

                     44.  Guido Weiss ( 1928-2021) 114;8722

                   45.  Bertram Kostant ( 1928-2017) 114;4367

                 46.  Andrei V. Zelevinsky ( 1953-2013) 113;8482

                 47.  Wilhelmus Luxemburg ( 1929-2018) 111;1107

                   48.  Nikolai Chernov ( 1956-2014) 111;1859

                  49.  Manfredo P. doCarmo ( 1928-2018) 110;3685

                    50.  Peter M. Gruber ( 1941-2017) 107;1867

                   51.  Thomas M. Liggett ( 1944-2020) 106;3918

                   52.  Ronald K. Getoor ( 1929-2017) 106;1850

                  53.  Richard V. Kadison ( 1925-2018) 105;3633

                54.  Elliott Ward Cheney, Jr. ( 1929-2016) 105;1642

                   55.  Ernest A. Michael ( 1925-2013) 102;2845

                     56.  Vera Pless ( 1931-2020) 101;2248

                    57.  Victor L. Shapiro ( 1924-2013) 99;239

                     58.  Herman Rubin ( 1926-2018) 99;498

                   59.  Isadore M. Singer ( 1924-2021) 96;5215

                      60.  John T. Tate ( 1925-2019) 95;5626

                    61.  Maurice H. Heins ( 1915-2015) 93;401

                     62.  Edward Odell ( 1947-2013) 92;1337

                    63.  Bjarni Jonsson ( 1920-2016) 90;1632

                    64.  Junichi Igusa ( 1924-2013) 89;1610

                    65.  Mary Ellen Rudin ( 1924-2013) 87;791

                     66.  Tom M. Apostol ( 1923-2016) 87;2886

                       67.  Lee Lorch ( 1915-2014) 86;272

                 68.  Clifford John Earle, Jr. ( 1935-2017) 84;1507

                     69.  Walter Craig ( 1953-2019) 84;2084

                70.  Cathleen Synge Morawetz ( 1923-2017) 82;1295

                    71.  Edward G. Effros ( 1935-2019) 82;3484

                   72.  Richard Pollack ( 1935-2018) 81;1658

                   73.  Heini Halberstam ( 1926-2014) 81;1086

                    74.  Charles J. Stone ( 1936-2019) 79;2604

                    75.  Joseph F. Traub ( 1932-2015) 78;1509

                   76.  Colin J. Bushnell ( 1947-2021) 77;1624

              77.  Komaravolu Chandrasekharan ( 1920-2017) 75;778

               78.  Jacques Jean-Pierre Neveu ( 1932-2016) 74;1172

                    79.  Robert R. Phelps ( 1926-2013) 71;2519

                   80.  Pierre E. Conner, Jr. ( 1932-2018) 71;1118

                 81.  Lloyd Stowell Shapley ( 1923-2016) 70;4107

                   82.  Elwyn Berlekamp ( 1940-2019) 70;1550

          83.  Conjeeveram Srirangachari Seshadri ( 1932-2020) 69;1726

                     84.  Samuel Gitler ( 1933-2014) 67;560

                    85.  William A. Veech ( 1938-2016) 66;1891

                    86.  Rudolf Kalman ( 1930-2016) 66;1826

                  87.  Jane Cronin Scanlon ( 1922-2018) 66;206

                      88.  Alan Baker ( 1939-2018) 65;1596

                    89.  Edward Nelson ( 1932-2014) 64;2035

                   90.  Albert Nijenhuis ( 1926-2015) 61;1054

                   91.  Joyce R. McLaughlin ( 1939-2017) 61;914

                     92.  Walter Noll ( 1925-2017) 59;1601

                    93.  Edward Fadell ( 1926-2018) 59;1065

                     94.  EdgarH Brown ( 1926-2021) 58;586

                    95.  Isaac Namioka ( 1928-2019) 57;1183

                     96.  Tim D. Cochran ( 1955-2014) 57;1265

                  97.  Donald L. Burkholder ( 1927-2013) 57;2165

                    98.  Steven M. Zucker ( 1949-2019) 55;756

                    99.  Robert T. Seeley ( 1932-2016) 53;1510

                  100.  Oscar E. Lanford, III ( 1940-2013) 50;734

                    101.  Robert M. Miura ( 1938-2018) 49;1370

                   102.  Richard L. Bishop ( 1931-2019) 49;1173

                      103.  John Roe ( 1959-2018) 48;2116

                 104.  Aldridge Bousfield ( 1941-2020) 48;2563

                105.  Clarence W. Wilkerson, Jr. ( 1944-2019) 47;749

                      106.  Roy Adler ( 1931-2016) 47;1231

                   107.  Richard D. Schafer ( 1918-2014) 46;712

                      108.  Ray Kunze ( 1928-2014) 46;452

                      109.  David Goss ( 1952-2017) 46;882

                   110.  Harold M. Edwards ( 1936-2020) 46;542

                     111.  David Cantor ( 1935-2012) 45;451

                     112.  Louis Howard ( 1929-2015) 42;275

                    113.  Paul Chernoff ( 1942-2017) 42;658

                     114.  Robert Ellis ( 1926-2013) 41;978

                    115.  Lynne H. Walling ( 1958-2021) 40;101

                     116.  Maurice Sion ( 1927-2018) 39;525

                  117.  Stephen A. Mitchell ( 1951-2017) 39;366

                 118.  Aderemi Oluyomi Kuku ( 1941-2022) 39;121

                    119.  Paul J. Sally, Jr. ( 1933-2013) 37;281

                     120.  John C. Moore ( 1923-2016) 37;1254

                   121.  Steve Armentrout ( 1930-2020) 37;145

                     122.  Howard Osborn ( 1928-2019) 36;52

                     123.  William Bade ( 1924-2012) 35;426

                     124.  Henry Wente ( 1936-2020) 34;526

                    125.  Kenneth I. Gross ( 1938-2017) 34;181

                    126.  Charles C. Sims ( 1937-2017) 33;503

                     127.  John Hempel ( 1935-2022) 31;1211

                    128.  Donald Babbitt ( 1936-2020) 31;161

                    129.  Kenneth I. Appel ( 1932-2013) 31;873

                    130.  Lesley M. Sibner ( 1934-2013) 29;157

                    131.  John F. Nash, Jr. ( 1928-2015) 29;3185

                  132.  Solomon Feferman ( 1928-2016) 29;1146

                     133.  John H. Walter ( 1927-2021) 27;214

                  134.  Alphonse T. Vasquez ( 1938-2021) 27;133

                  135.  J. Michael Boardman ( 1938-2021) 26;781

                  136.  Samuel M. Rankin, III ( 1945-2020) 25;145

                    137.  Thomas Nevins ( 1972-2020) 24;173

                    138.  Myles Tierney ( 1937-2017) 23;428

                   139.  Robert J. Blattner ( 1931-2015) 23;336

                    140.  James M. Kister ( 1930-2018) 19;142

                141.  Arthur Herbert Copeland ( 1926-2019) 19;41

                   142.  Lynn Arthur Steen ( 1941-2015) 18;341

                   143.  Arthur P. Mattuck ( 1930-2021) 17;152

                  144.  Carole B. Lacampagne ( 1933-2021) 11;12

                      145.  Lida K. Barrett ( 1927-2021) 5;2

{\bf Appendix III: Ranking the  AMS Fellows who died before Aug. 2022 According to Decreasing Number of Citations}

                    1.  Elias M. Stein ( 1931-2018) 237;27492

                   2.  Louis Nirenberg ( 1925-2020) 175;17621

                   3.  Lars Hormander ( 1931-2012) 159;15121

                  4.  Sir Michael Atiyah ( 1929-2019) 269;12281

                    5.  Andrew J. Majda ( 1949-2021) 323;11987

               6.  Jonathan Michael Borwein ( 1951-2016) 453;9810

                     7.  Guido Weiss ( 1928-2021) 114;8722

                  8.  Andrei V. Zelevinsky ( 1953-2013) 113;8482

                    9.  Ronald Graham ( 1935-2020) 363;7500

                   10.  Hans Weinberger ( 1928-2017) 115;6664

                    11.  Anatole Katok ( 1944-2018) 122;6030

                    12.  Jim Douglas, Jr. ( 1927-2016) 180;5937

                      13.  John T. Tate ( 1925-2019) 95;5626

                   14.  Isadore M. Singer ( 1924-2021) 96;5215

                    15.  Richard Askey ( 1933-2019) 163;4802

                    16.  Harry Kesten ( 1931-2016) 193;4772

                   17.  Felix E. Browder ( 1927-2016) 200;4743

                  18.  Robert Strichartz ( 1943-2021) 189;4633

                     19.  Peter Duren ( 1935-2020) 126;4490

                   20.  Bertram Kostant ( 1928-2017) 114;4367

                 21.  Lloyd Stowell Shapley ( 1923-2016) 70;4107

                    22.  Joel A. Smoller ( 1935-2017) 181;4096

                    23.  George R. Sell ( 1937-2015) 132;4068

                    24.  Harold Widom ( 1932-2021) 167;4066

                   25.  Vaughan F. R. Jones ( 1952-2020) 117;3924

                   26.  Thomas M. Liggett ( 1944-2020) 106;3918

                   27.  Branko Grunbaum ( 1929-2018) 271;3814

                   28.  Ronald G. Douglas ( 1938-2018) 168;3732

                  29.  Manfredo P. doCarmo ( 1928-2018) 110;3685

                   30.  Joseph B. Keller ( 1923-2016) 321;3655

                  31.  Richard V. Kadison ( 1925-2018) 105;3633

                      32.  Paul Fife ( 1930-2014) 118;3509

                    33.  Edward G. Effros ( 1935-2019) 82;3484

                     34.  John F. Nash, Jr. ( 1928-2015) 29;3185

                    35.  Robin Thomas ( 1962-2020) 122;2898

                     36.  Tom M. Apostol ( 1923-2016) 87;2886

                   37.  Eugene B. Dynkin ( 1924-2014) 222;2851

                   38.  Ernest A. Michael ( 1925-2013) 102;2845

                   39.  Richard M. Dudley ( 1938-2020) 138;2776

                    40.  Jan C. Willems ( 1939-2013) 186;2677

                 41.  Jean-Pierre Kahane ( 1926-2017) 288;2639

                    42.  Charles J. Stone ( 1936-2019) 79;2604

                  43.  Aldridge Bousfield ( 1941-2020) 48;2563

                    44.  Robert R. Phelps ( 1926-2013) 71;2519

               45.  Veeravalli S. Varadarajan ( 1937-2019) 149;2336

                     46.  Vera Pless ( 1931-2020) 101;2248

                  47.  Donald L. Burkholder ( 1927-2013) 57;2165

                   48.  Lawrence Shepp ( 1936-2013) 175;2137

                  49.  Murray Rosenblatt ( 1926-2019) 178;2128

                       50.  John Roe ( 1959-2018) 48;2116

                   51.  Andras Hajnal ( 1931-2016) 164;2116

                 52.  Shreeram Abhyankar ( 1930-2012) 190;2111

                     53.  Walter Craig ( 1953-2019) 84;2084

                    54.  Edward Nelson ( 1932-2014) 64;2035

                  55.  Gilbert Baumslag ( 1933-2014) 166;1937

                    56.  William A. Veech ( 1938-2016) 66;1891

                    57.  Peter M. Gruber ( 1941-2017) 107;1867

                   58.  Nikolai Chernov ( 1956-2014) 111;1859

                   59.  Ronald K. Getoor ( 1929-2017) 106;1850

                    60.  Rudolf Kalman ( 1930-2016) 66;1826

                  61.  Frederick R. Cohen ( 1946-2022) 136;1781

                   62.  W. Wistar Comfort ( 1933-2016) 153;1746

          63.  Conjeeveram Srirangachari Seshadri ( 1932-2020) 69;1726

                   64.  Georgia Benkart ( 1949-2022) 133;1692

                   65.  Richard Pollack ( 1935-2018) 81;1658

                66.  Elliott Ward Cheney, Jr. ( 1929-2016) 105;1642

                    67.  Bjarni Jonsson ( 1920-2016) 90;1632

                   68.  Colin J. Bushnell ( 1947-2021) 77;1624

                    69.  Junichi Igusa ( 1924-2013) 89;1610

                     70.  Walter Noll ( 1925-2017) 59;1601

                      71.  Alan Baker ( 1939-2018) 65;1596

                   72.  Elwyn Berlekamp ( 1940-2019) 70;1550

                    73.  Robert T. Seeley ( 1932-2016) 53;1510

                    74.  Joseph F. Traub ( 1932-2015) 78;1509

                 75.  Clifford John Earle, Jr. ( 1935-2017) 84;1507

                   76.  Solomon Golomb ( 1932-2016) 146;1497

                    77.  Robert M. Miura ( 1938-2018) 49;1370

                     78.  Edward Odell ( 1947-2013) 92;1337

                 79.  Vladimir F. Demyanov ( 1938-2014) 183;1309

                80.  Cathleen Synge Morawetz ( 1923-2017) 82;1295

                    81.  Mark Mahowald ( 1931-2013) 164;1281

                     82.  Tim D. Cochran ( 1955-2014) 57;1265

                     83.  John C. Moore ( 1923-2016) 37;1254

                      84.  Roy Adler ( 1931-2016) 47;1231

                     85.  John Hempel ( 1935-2022) 31;1211

                    86.  Isaac Namioka ( 1928-2019) 57;1183

                   87.  Richard L. Bishop ( 1931-2019) 49;1173

               88.  Jacques Jean-Pierre Neveu ( 1932-2016) 74;1172

                   89.  Solomon Feferman ( 1928-2016) 29;1146

                   90.  Sibe Mardesic ( 1927-2016) 159;1130

                   91.  Pierre E. Conner, jr. ( 1932-2018) 71;1118

                 92.  Wilhelmus Luxemburg ( 1929-2018) 111;1107

                   93.  Heini Halberstam ( 1926-2014) 81;1086

                    94.  Edward Fadell ( 1926-2018) 59;1065

                   95.  Albert Nijenhuis ( 1926-2015) 61;1054

                     96.  Robert Ellis ( 1926-2013) 41;978

                   97.  Joyce R. McLaughlin ( 1939-2017) 61;914

                      98.  David Goss ( 1952-2017) 46;882

                    99.  Kenneth I. Appel ( 1932-2013) 31;873

                   100.  Mary Ellen Rudin ( 1924-2013) 87;791

                  101.  J. Michael Boardman ( 1938-2021) 26;781

              102.  Komaravolu Chandrasekharan ( 1920-2017) 75;778

                    103.  Steven M. Zucker ( 1949-2019) 55;756

                104.  Clarence W. Wilkerson, Jr. ( 1944-2019) 47;749

                  105.  Oscar E. Lanford, III ( 1940-2013) 50;734

                   106.  Richard D. Schafer ( 1918-2014) 46;712

                     107.  Harvey Cohn ( 1923-2014) 118;664

                    108.  Paul Chernoff ( 1942-2017) 42;658

                     109.  Edgar H. Brown ( 1926-2021) 58;586

                    110.  Samuel Gitler ( 1933-2014) 67;560

                   111.  Harold M. Edwards ( 1936-2020) 46;542

                     112.  Henry Wente ( 1936-2020) 34;526

                     113.  Maurice Sion ( 1927-2018) 39;525

                    114.  Charles C. Sims ( 1937-2017) 33;503

                     115.  Herman Rubin ( 1926-2018) 99;498

                      116.  Ray Kunze ( 1928-2014) 46;452

                     117.  David Cantor ( 1935-2012) 45;451

                    118.  Myles Tierney ( 1937-2017) 23;428

                     119.  William Bade ( 1924-2012) 35;426

                    120.  Maurice H. Heins ( 1915-2015) 93;401

                   121.  Donald A. S. Fraser ( 1925-2020) 182;400

                  122.  Stephen A. Mitchell ( 1951-2017) 39;366

                   123.  Lynn Arthur Steen ( 1941-2015) 18;341

                   124.  Robert J. Blattner ( 1931-2015) 23;336

                    125.  Paul J. Sally, Jr. ( 1933-2013) 37;281

                     126.  Louis Howard ( 1929-2015) 42;275

                      127.  Lee Lorch ( 1915-2014) 86;272

                   128.  Victor L. Shapiro ( 1924-2013) 99;239

                     129.  John H. Walter ( 1927-2021) 27;214

                  130.  Jane Cronin Scanlon ( 1922-2018) 66;206

                    131.  Kenneth I. Gross ( 1938-2017) 34;181

                    132.  Thomas Nevins ( 1972-2020) 24;173

                    133.  Donald Babbitt ( 1936-2020) 31;161

                    134.  Lesley M. Sibner ( 1934-2013) 29;157

                   135.  Arthur P. Mattuck ( 1930-2021) 17;152

                  136.  Samuel M. Rankin, III ( 1945-2020) 25;145

                   137.  Steve Armentrout ( 1930-2020) 37;145

                    138.  James M. Kister ( 1930-2018) 19;142

                  139.  Alphonse T. Vasquez ( 1938-2021) 27;133

                 140.  Aderemi Oluyomi Kuku ( 1941-2022) 39;121

                    141.  Lynne H. Walling ( 1958-2021) 40;101

                     142.  Howard Osborn ( 1928-2019) 36;52

                143.  Arthur Herbert Copeland ( 1926-2019) 19;41

                  144.  Carole B. Lacampagne ( 1933-2021) 11;12

                      145.  Lida K. Barrett ( 1927-2021) 5;2

{\bf Acknowledgment}: Many thanks are due to Shalosh B. Ekhad for its dilligent computations, that were generated using the Maple package

{\tt https://sites.math.rutgers.edu/\~{}zeilberg/tokhniot/AMSrip.txt} \quad .

{\bf References}

[1] The AMS, List of AMS Fellows {\tt http://www.ams.org/cgi-bin/fellows/fellows.cgi} (Viewed Aug. 31, 2022)

[2] MathSciNet, {\tt https://mathscinet.ams.org/mathscinet/} (requires subscription)

\bigskip
\hrule
\bigskip
Corresponding author:  Doron Zeilberger, Department of Mathematics, Rutgers University (New Brunswick), Hill Center-Busch Campus, 110 Frelinghuysen
Rd., Piscataway, NJ 08854-8019, USA. \hfill\break
Email: {\tt DoronZeil] at gmail dot com}   \quad .

Written: {\bf Sept. 29, 2022}. 


\end